\documentclass[12pt,twoside,reqno]{amsart}
\usepackage{amsmath}
\usepackage{amsfonts}
\usepackage{amssymb}
\usepackage{color}
\usepackage{mathrsfs}
\newcommand{\beqn}{\begin{equation}}
\newcommand{\eeqn}{\end{equation}}
\newcommand{\bean}{\begin{eqnarray}}
\newcommand{\eean}{\end{eqnarray}}
\DeclareMathAlphabet{\mathpzc}{OT1}{pzc}{m}{it}
\newtheorem{theorem}{Theorem}

\newtheorem{lemma}[theorem]{Lemma}
\newtheorem{proposition}[theorem]{Proposition}

\newtheorem{remark}[theorem]{\it Remark}
\begin{document}

\title[Finite time singularity in a free boundary problem modeling MEMS]{Finite time singularity in a free boundary problem\\ modeling MEMS}

\author{Joachim Escher}
\address{Leibniz Universit\"at Hannover\\ Institut f\" ur Angewandte Mathematik \\ Welfengarten 1 \\ D--30167 Hannover\\ Germany}
\email{escher@ifam.uni-hannover.de}

\author{Philippe Lauren\c{c}ot}\thanks{Partially supported by the Deutscher Akademischer Austausch Dienst (DAAD)}
\address{Institut de Math\'ematiques de Toulouse, CNRS UMR~5219, Universit\'e de Toulouse \\ F--31062 Toulouse Cedex 9, France}
\email{laurenco@math.univ-toulouse.fr}

\author{Christoph Walker}
\address{Leibniz Universit\"at Hannover\\ Institut f\" ur Angewandte Mathematik \\ Welfengarten 1 \\ D--30167 Hannover\\ Germany}
\email{walker@ifam.uni-hannover.de}

\date{\today}

\begin{abstract}
The occurrence of a finite time singularity is shown for a free boundary problem modeling microelectromechanical systems (MEMS) when the applied voltage exceeds some value. The model involves a singular nonlocal reaction term and a nonlinear curvature term accounting for large deformations.
\end{abstract}

\maketitle

\section{Introduction}
\label{s.int}

An idealized electostatically actuated microelectromechanical system (MEMS) consists of a fixed horizontal ground plate held at zero potential above which an elastic membrane held at potential $V$ is suspended. A Coulomb force is generated by the potential difference across the device and results in a deformation of the membrane, thereby converting  electrostatic energy into mechanical energy, see \cite{BP11,EGG10,PB03} for a more detailed account and further references. After a suitable scaling and assuming homogeneity in transversal horizontal direction, the ground plate is assumed to be located at $z=-1$ and the membrane displacement $u=u(t,x)\in (-1,\infty)$ with $t>0$ and $x\in I:=(-1,1)$ evolves according to
\begin{equation}\label{u}
\partial_t u - \partial_x\left(\frac{\partial_x u}{\sqrt{1+\varepsilon^2(\partial_x u)^2}}\right) = - \lambda \left( \varepsilon^2\ |\partial_x\psi(t,x,u(t,x))|^2 + |\partial_z\psi(t,x,u(t,x))|^2 \right)\ ,
\end{equation}
for $t>0$ and $x\in I$
with boundary conditions
\begin{equation}\label{bcu}
u(t,\pm 1)=0\ ,\quad t>0\ ,
\end{equation}
and initial condition
\begin{equation}\label{ic}
u(0,x)=u^0(x)\ ,\quad x\in I\ .
\end{equation}
The electrostatic potential $\psi=\psi(t,x,z)$ satisfies a rescaled Laplace equation in the region 
$$
\Omega(u(t)) := \left\{ (x,z)\in I\times (-1,\infty)\ :\ -1 < z < u(t,x) \right\}
$$ 
between the plate and the membrane which reads
\begin{align}
\varepsilon^2 \,\partial_x^2\psi+\partial_z^2\psi & = 0 \ ,\quad (x,z)\in \Omega(u(t))\ ,\quad t>0\ , \label{psi} \\
\psi(t,x,z) & = \frac{1+z}{1+u(t,x)}\ ,\quad (x,z)\in\partial\Omega(u(t))\ ,\quad t>0\ ,\label{bcpsi}
\end{align}
where $\varepsilon>0$ denotes the aspect ratio of the device and $\lambda>0$ is proportional to the square of the applied voltage. The dynamics of $(u,\psi)$ is thus given by the coupling of a quasilinear parabolic equation for $u$ and an elliptic equation in a moving domain for $\psi$, the latter being only well-defined as long as the membrane does not touch down on the ground plate, that is, $u$ does not reach the value $-1$. To guarantee optimal operating conditions of the device, this touchdown phenomenon has to be controlled and its occurrence is obviously related to the value of $\lambda$. 

The main difficulty to be overcome in the analysis of \eqref{u}-\eqref{bcpsi} is the nonlocal and nonlinear implicit dependence on $u$ of the right-hand side of \eqref{u} which is also singular if $u$ approaches $-1$. Except for the singularity, these features disappear when setting $\varepsilon=0$ in \eqref{u}-\eqref{bcpsi}, a commonly made assumption which reduces \eqref{u}-\eqref{bcpsi}  to a singular semilinear reaction-diffusion equation. This so-called small aspect ratio model has received considerable attention in recent years, see \cite{EGG10,PB03} and the references therein. In this simplified situation, it has been established that touchdown does not take place if $\lambda$ is below a certain threshold value $\lambda_*>0$, but occurs if $\lambda$ exceeds this value \cite{EGG10,FMPS07,GPW05}. 

We have recently investigated the well-posedness of \eqref{u}-\eqref{bcpsi} and established the following result~\cite{ELWyy}.

\medskip

\begin{theorem}[{\bf Local Well-Posedness}]\label{A}
Let $q\in (2,\infty)$, $\varepsilon>0$, $\lambda>0$, and consider an initial value 
\begin{equation}
u^0\in  W_q^2(I) \;\;\text{ such that }\;\; u^0(\pm1)=0 \;\;\text{ and }\;\; 0 \ge u^0(x)>-1 \;\;\text{ for }\;\; x\in I\ . \label{aic}
\end{equation} 
Then there is a unique maximal solution $(u,\psi)$ to \eqref{u}-\eqref{bcpsi} on the maximal interval of existence $[0,T_m^\varepsilon)$ in the sense that
$$
u\in C^1\big([0,T_m^\varepsilon),L_q(I)\big)\cap C\big([0,T_m^\varepsilon), W_q^2(I)\big)
$$
satisfies \eqref{u}-\eqref{ic} together with
\begin{equation}
0 \ge u(t,x)>-1\ ,\quad (t,x)\in [0,T_m^\varepsilon)\times I\ , \label{b0}
\end{equation}
and $\psi(t)\in W^2_{2}\big(\Omega(u(t))\big)$ solves \eqref{psi}-\eqref{bcpsi} for each $t\in [0,T_m^\varepsilon)$.
\end{theorem}

\medskip

We have also shown in \cite{ELWyy} that, if $\lambda$ and $u^0$ are sufficiently small, the solution $(u,\psi)$ to \eqref{u}-\eqref{bcpsi} exists for all times (i.e. $T_m^\varepsilon=\infty$) and touchdown does not take place, not even in infinite time.

\medskip

\begin{theorem}[{\bf Global Existence}]\label{B} 
Let $q\in (2,\infty)$, $\varepsilon>0$, and consider an initial value $u^0$ satisfying~\eqref{aic}. Given $\kappa\in (0,1)$, there are $\lambda_*(\kappa)>0$ and $r(\kappa)>0$ such that, if $\lambda\in (0,\lambda_*(\kappa))$ and $\|u^0\|_{W_q^2(I)}\le r(\kappa)$, the maximal solution $(u,\psi)$ to \eqref{u}-\eqref{bcpsi} exists for all times and $u(t,x)\ge -1+\kappa$ for $(t,x)\in [0,\infty)\times I$.
\end{theorem}

\medskip

On the other hand, we have been able to prove that no stationary solution to \eqref{u}-\eqref{bcpsi} exists provided $\lambda$ is sufficiently large. However,whether or not $T_m^\varepsilon$ is finite in this case has been left as an open question. The purpose of this note is to show that -- as expected on physical grounds -- $T_m^\varepsilon$ is indeed finite for $\lambda$ sufficiently large.

\medskip

\begin{theorem}[{\bf Finite time singularity}]\label{C}
Let $q\in (2,\infty)$, $\varepsilon>0$, and consider an initial value $u^0$ satisfying~\eqref{aic}. If $\lambda>1/\varepsilon$ and $(u,\psi)$ denotes the maximal solution to \eqref{u}-\eqref{bcpsi} defined on $[0,T_m^\varepsilon)$, then $T_m^\varepsilon<\infty$.
\end{theorem}

\medskip

The criterion $\lambda>1/\varepsilon$ is likely to be far from optimal. As we shall see below, improving it would require to have a better control on $\partial_x u(\pm 1)$. The proof of Theorem~\ref{C} relies on the derivation of a chain of estimates which allow us to obtain a lower bound on the $L_1$-norm of the right-hand side of \eqref{u} depending only on $u$. The lower bound thus obtained is in fact the mean value of a convex function of $u$, and we may then end the proof with the help of  Jensen's inequality, an argument which has already been used for the small aspect ratio model, see \cite{FMPS07,GPW05}.

\medskip

We shall point out that, in contrast to the small aspect ratio model, the finiteness of $T_m^\varepsilon$ does not guarantee that the touchdown phenomenon really takes place as $t\to T_m^\varepsilon$.  Indeed, according to \cite[Theorem~1.1~(ii)]{ELWyy}, the finiteness of $T_m^\varepsilon$ implies that $\min_{[-1,1]}{u(t)}\longrightarrow -1$ or $\|u(t)\|_{W_q^2(I)} \longrightarrow \infty$ as $t\to T_m^\varepsilon$. While the former corresponds to the touchdown behaviour, the latter is more likely to be interpreted as the membrane being no longer the graph of a function at time $T_m^\varepsilon$.

\section{Proof of Theorem~\ref{C}}
\label{s.proof}

Let $q\in (2,\infty)$, $\varepsilon>0$, $\lambda>0$ and consider an initial value $u^0$ satisfying \eqref{aic}. We denote the maximal solution to \eqref{u}-\eqref{bcpsi} defined on $[0,T_m^\varepsilon)$ by $(u,\psi)$. Differentiating the boundary conditions \eqref{bcpsi}, we readily obtain 
\begin{equation}
\partial_x \psi(t,x,-1) = \partial_x \psi(t,x,u(t,x)) + \partial_x u(t,x)\, \partial_z \psi(t,x,u(t,x)) = 0\ , \qquad (t,x)\in (0,T_m^\varepsilon)\times I\ , \label{b3a} 
\end{equation}
and
\begin{equation}
\partial_z \psi(t,\pm 1,z) = 1\ ,  \qquad (t,z)\in (0,T_m^\varepsilon)\times (-1,0)\ . \label{b3b}
\end{equation}

Additional information on the boundary behaviour of $\psi$ is provided by the next lemma.

\medskip

\begin{lemma}\label{le0}
For $t\in (0,T_m^\varepsilon)$,
\begin{align}
1+z \le \psi(t,x,z) \le 1\ , & \qquad (x,z)\in \Omega(u(t))\ , \label{b1} \\
\pm \partial_x\psi(t,\pm 1, z) \le 0\ , & \qquad z\in (-1,0)\ . \label{b2}
\end{align}
\end{lemma}

\noindent\textbf{Proof.}
Fix $t\in (0,T_m^\varepsilon)$. The upper bound in \eqref{b1} readily follows from the maximum principle. Next, the function $\sigma$, defined by $\sigma(x,z)=1+z$, obviously satisfies $\varepsilon^2 \partial_x^2 \sigma + \partial_z^2 \sigma = 0$ in $\Omega(u(t))$ as well as 
\begin{align*}
\sigma(\pm 1,z) & = 1+z = \psi(t,\pm 1, z)\ , & \qquad z\in (-1,0)\ , \\
\sigma(x,-1) & = 0 = \psi(t,x,-1)\ , & \qquad x\in (-1,1)\ .
\end{align*}
Owing to the non-positivity \eqref{b0} of $u(t)$, it also satisfies
$$
\sigma(x,u(t,x)) = 1+u(t,x) \le 1 = \psi(t,x,u(t,x))\ , \qquad x\in (-1,1)\ , 
$$
and we infer from the comparison principle that $\psi(t,x,z)\ge \sigma(x,z)$ for $(x,z)\in\Omega(u(t))$. It then follows from \eqref{b1} that $\psi(t,x,z)\ge 1+z = \psi(t,\pm 1,z)$ for $(x,z)\in\Omega(u(t))$ which readily implies \eqref{b2}. \hfill $\square$

\medskip

To simplify notations, we set
\begin{equation}
\gamma_m(t,x) := \partial_z \psi(t,x,u(t,x))\ , \quad \gamma_g(t,x) := \partial_z \psi(t,x,-1)\ , \qquad (t,x)\in (0,T_m^\varepsilon)\times (-1,1)\ , \label{b4}
\end{equation}
and first derive an upper bound of the $L_1$-norm of the right-hand side of \eqref{u}, observing that, due to \eqref{b3a}, it also reads
$$
- \lambda \varepsilon^2\ |\partial_x\psi(t,x,u(t,x))|^2 + |\partial_z\psi(t,x,u(t,x))|^2 = - \lambda \left( 1 + \varepsilon^2 (\partial_x u(t,x))^2 \right) \gamma_m(t,x)^2\ .
$$

\begin{lemma}\label{le1} 
For $t\in (0,T_m^\varepsilon)$, 
\begin{equation}
\int_{-1}^1 \left( 1 + \varepsilon^2 (\partial_x u(t,x))^2 \right) \gamma_m(t,x)^2\ \mathrm{d}x \ge 2 \int_{-1}^1 \left( 1 + \varepsilon^2 (\partial_x u(t,x))^2 \right) \gamma_m(t,x)\ \mathrm{d}x - 2\ . \label{b5}
\end{equation}
\end{lemma}

\noindent\textbf{Proof.} Fix $t\in (0,T_m^\varepsilon)$. We multiply \eqref{psi} by $\partial_z\psi(t)-1$ and integrate over $\Omega(u(t))$. Using \eqref{b3a}, \eqref{b3b}, and Green's formula we obtain
\begin{align*}
0 = & - \varepsilon^2 \int_{\Omega(u)} \partial_x\partial_z \psi\ \partial_x \psi\ \mathrm{d}(x,z) + \varepsilon^2 \int_{-1}^1 (\partial_x u)^2\ \gamma_m \left( \gamma_m -1 \right)\ \mathrm{d}x \\
& - \frac{1}{2} \int_{-1}^1 \left( \gamma_g^2 - 2 \gamma_g \right)\ \mathrm{d}x + \frac{1}{2} \int_{-1}^1 \left( \gamma_m^2 - 2 \gamma_m \right)\ \mathrm{d}x \ .
\end{align*}
Since
$$
\int_{\Omega(u)} \partial_x\partial_z \psi\ \partial_x \psi\ \mathrm{d}(x,z) = \frac{1}{2} \int_{-1}^1 (\partial_x u)^2\ \gamma_m^2\ \mathrm{d}x 
$$
by \eqref{b3a} and since $\gamma_g^2 - 2 \gamma_g \ge -1$, we end up with \eqref{b5}. \hfill $\square$

\medskip

We again use \eqref{psi} to obtain a lower bound for the boundary integral of the right-hand side of \eqref{b5} which depends on the Dirichlet energy of $\psi$.

\medskip

\begin{lemma}\label{le2}
For $t\in (0,T_m^\varepsilon)$, 
\begin{equation}
\int_{-1}^1 \left( 1 + \varepsilon^2 (\partial_x u(t,x))^2 \right) \gamma_m(t,x)\ \mathrm{d}x \ge \int_{\Omega(u(t))} \left( \varepsilon^2 \vert \partial_x \psi(t,x,z)\vert^2 + \vert \partial_z \psi(t,x,z)\vert^2 \right)\ \mathrm{d}(x,z)\ . \label{b6}
\end{equation}
\end{lemma}

\noindent\textbf{Proof.} Fix $t\in (0,T_m^\varepsilon)$. We multiply \eqref{psi} by $\psi(t)$ and integrate over $\Omega(u(t))$. Using \eqref{bcpsi}, \eqref{b3a}, and Green's formula we obtain
\begin{align*}
0 = & - \int_{\Omega(u(t))} \left( \varepsilon^2 \vert \partial_x \psi(t,x,z)\vert^2 + \vert \partial_z \psi(t,x,z)\vert^2 \right)\ \mathrm{d}(x,z) + \varepsilon^2 \int_{-1}^0 (1+z)\ \partial_x \psi(t,1,z)\ \mathrm{d}z \\
& - \varepsilon^2 \int_{-1}^0 (1+z)\ \partial_x \psi(t,-1,z)\ \mathrm{d}z + \varepsilon^2 \int_{-1}^1 (\partial_x u(t,x))^2\ \gamma_m(t,x)\ \mathrm{d}x + \int_{-1}^1 \gamma_m(t,x)\ \mathrm{d}x \ .
\end{align*}
Owing to \eqref{b2}, the second and third terms of the right-hand side of the above equality are non-positive, whence \eqref{b6}. \hfill $\square$

\medskip

We finally argue as in \cite[Lemma~9]{ELWxx} to establish a connection between the Dirichlet energy of $\psi$ and $u$.

\medskip

\begin{lemma}\label{le3}
For $t\in (0,T_m^\varepsilon)$, 
\begin{equation}
\int_{\Omega(u(t))} \left( \varepsilon^2 \vert \partial_x \psi(t,x,z)\vert^2 + \vert \partial_z \psi(t,x,z)\vert^2 \right)\ \mathrm{d}(x,z) \ge \int_{-1}^1 \frac{\mathrm{d}x}{1+u(t,x)} \ . \label{b7}
\end{equation}
\end{lemma}

\noindent\textbf{Proof.}  Let $t\in (0,T_m^\varepsilon)$ and $x\in (-1,1)$. We deduce from \eqref{bcpsi} and the Cauchy-Schwarz inequality that
\begin{align}
\frac{1}{1+u(t,x)} = & \frac{\left( \psi(t,x,u(t,x)) - \psi(t,x,-1) \right)^2}{1+u(t,x)} = \frac{1}{1+u(t,x)} \left( \int_{-1}^{u(t,x)} \partial_z \psi(t,x,z)\ \mathrm{d}z \right)^2 \nonumber\\
\le & \int_{-1}^{u(t,x)} \left( \partial_z \psi(t,x,z) \right)^2\ \mathrm{d}z \ . \label{b7b}
\end{align}
Integrating the above inequality with respect to $x\in (-1,1)$ readily gives \eqref{b7}. \hfill $\square$

\medskip

\begin{remark}\label{re1} Observe that \eqref{b7b} provides a quantitative estimate on the singularity of $\partial_z\psi$ generated by $u$ when touchdown occurs.
\end{remark}

\medskip

Combining the three lemmas above with Jensen's inequality give the following estimate.

\medskip

\begin{proposition}\label{pr4}
For $t\in (0,T_m^\varepsilon)$, 
\begin{equation}
\int_{-1}^{1} \left( 1 + \varepsilon^2 (\partial_x u(t,x))^2 \right) \gamma_m(t,x)^2\ \mathrm{d}x \ge 4 \varphi\left( \frac{1}{2} \int_{-1}^1 u(t,x)\ \mathrm{d}x \right) - 2\ , \label{b8}
\end{equation}
where $\varphi(r):=1/(1+r)$, $r\in (-1,\infty)$.
\end{proposition}

\noindent\textbf{Proof.} Fix $t\in (0,T_m^\varepsilon)$. We infer from Lemma~\ref{le1}, Lemma~\ref{le2}, and Lemma~\ref{le3} that
$$
\int_{-1}^{1} \left( 1 + \varepsilon^2 (\partial_x u(t,x))^2 \right) \gamma_m(t,x)^2\ \mathrm{d}x \ge 2 \int_{-1}^1 \varphi(u(t,x))\ \mathrm{d}x - 2\ .
$$
To complete the proof, we argue as in \cite{FMPS07,GPW05} and use the convexity of $\varphi$ and Jensen's inequality to obtain~\eqref{b8}. \hfill $\square$

\bigskip

\noindent\textbf{Proof of Theorem~\ref{C}.} Introducing 
$$
E(t) := - \frac{1}{2} \int_{-1}^1 u(t,x)\ \mathrm{d}x\ , \qquad t\in [0,T_m^\varepsilon)\ ,
$$
the bounds \eqref{b0} ensure that
\begin{equation}
0 \le E(t) < 1\ , \qquad t\in [0,T_m^\varepsilon)\ . \label{b9}
\end{equation}
It follows from \eqref{u}, \eqref{b3a},  and Proposition~\ref{pr4} that
\begin{align}
\frac{dE}{dt}(t) = & - \frac{1}{2} \left[ \frac{\partial_x u(t,x)}{\sqrt{1 + \varepsilon^2 \left( \partial_x u(t,x) \right)^2}} \right]_{x=-1}^{x=1} + \frac{\lambda}{2} \int_{-1}^1 \left( 1 + \varepsilon^2 (\partial_x u(t,x))^2 \right) \gamma_m(t,x)^2\ \mathrm{d}x \nonumber \\
\ge & \ F_\lambda(E) := 2\lambda \varphi(-E) - \lambda - \frac{1}{\varepsilon}\ . \label{b10}
\end{align}
If $\lambda>1/\varepsilon$, we note that $F_\lambda(0)>0$ and thus $F_\lambda(r)\ge F_\lambda(0)>0$ for $r\in [0,1)$ due to the monotonicity of $F_\lambda$. Since $E(0)\ge 0$ by \eqref{b9}, it follows from \eqref{b10} and the properties of $F_\lambda$ that $t\mapsto E(t)$ is increasing on $[0,T_m^\varepsilon)$. Consequently,
$$
\frac{dE}{dt}(t) \ge F_\lambda(E(0)) \ge F_\lambda(0)\ , \qquad t\in [0,T_m^\varepsilon)\ .
$$
Integrating the previous inequality with respect to time and using \eqref{b9}, we end up with the inequality $1\ge E(0) + F_\lambda(0) T_m^\varepsilon$ which provides the claimed finiteness of $T_m^\varepsilon$. \hfill $\square$
 
\section*{Acknowledgements}

This research was done while Ph.L. was enjoying the kind hospitality of the Institut f\"ur Angewandte Mathematik of the Leibniz Universit\"at Hannover.



\end{document}